\documentclass{article}
\usepackage{graphicx}
\graphicspath{ {./graphics-bw/} }
\begin{document}

\title{Associative submanifolds and gradient cycles}
\author{Simon Donaldson and Christopher Scaduto}
\date{\today}
\maketitle

\newcommand{\bR}{{\bf R}}
\newcommand{\bC}{{\bf C}}
\newtheorem{cor}{Corollary}
  \ \ \ \ \ \ \ \ \  {\it Dedicated to Professor S-T.  Yau, for his 70th birthday}
  
  \
  
  \

\section{Introduction}
  There is a standard cross product $\bR^{7}\times \bR^{7}\rightarrow \bR^{7}$
which is related to the $8$-dimensional Cayley algebra in the same way as
the familiar cross product on $\bR^{3}$ is related to the quaternion algebra.
The Euclidean form on $\bR^{7}$ can be recovered from the cross product by
taking $-1/6$ times  $ {\rm Tr}\left( L_{x}\right)^{2}$ where
  $L_{x}(y)= x\times y$. 
    One definition of a  $G_{2}$ structure on a $7$-manifold $M$ is a cross
product $TM\times TM\rightarrow TM$ which is equivalent to the standard model
at each point. The $G_{2}$ structure is torsion-free if the product is parallel
with respect to the Levi-Civita connection of the Riemannian metric on $M$
induced by the cross product, as above. A $3$-dimensional submanifold $P\subset
M$ is called {\it associative} if its tangent spaces are closed under the
cross-product.  Associative submanifolds are interesting from many points
of view. They are examples of calibrated submanifolds in the sense of Harvey and Lawson \cite{kn:HL} and they are fundamental objects in 
$G_{2}$-geometry. 
The purpose of this article is to  explore a model for associative
submanifolds in a class of $G_{2}$-structures near to an \lq\lq adiabatic
limit'', which was discussed in \cite{kn:D1}. (Related ideas were introduced earlier by Gukov, Yau and Zaslow in \cite{kn:GYZ}.) In this case the $7$-manifold
$M$ is equipped with a fibration $\pi:M \rightarrow N$  (with some singular
fibres) over a $3$-dimensional base and the smooth fibres are diffeomorphic
to  K3 surfaces.  Our model is built on the precise knowledge of complex curves
in the fibres coming from the standard theory of K3 surfaces. The basic idea
is to describe associative submanifolds via certain graphs in the base $N$.
This idea is in the same vein to other constructions in the literature, in
particular of tropical curves in Calabi-Yau manifolds with Strominger-Yau-Zaslow
fibrations (see the futher discussion in Section 5 below).
    
    We should emphasise at the outset that in this article we only take the
first steps towards a more comprehensive theory that one can hope will emerge
in the future. In particular we do not prove anything about actual associative
submanifolds here. The main purpose of this article is to develop an independent
 \lq\lq adiabatic limit theory'', involving graphs in $3$-manifolds, and
to show that it can mimic important known phenomena of singularity formation
for associative submanifolds\index{}.

This work was supported by the Simons Foundation through the Simons Collaboration {\it Special holonomy in geometry, analysis and physics}. The authors are grateful to Rodrigo Barbosa, Chris Gerig and Andrew Neitzke for helpful discussions.

    \section{Review of standard theory and Joyce's conjecture}
    
    We will now review more systematically some standard material on $G_{2}$-geometry.

    A $G_{2}$ structure on  $M$ defines a $3$-form $\phi$, related to the
cross-product and metric by
    $$   \phi(\xi_{1}, \xi_{2},\xi_{3})= \langle \xi_{1}\times \xi_{2}, \xi_{3}\rangle.
$$
    We also have the Hodge dual $4$-form $*\phi$. If the structure is torsion-free these forms are parallel and hence closed.  An alternative formulation
of the associative condition for a $3$-dimensional submanifold $P$ is that
for each point $x$ of $P$ and tangent vector $v\in TM_{x}$ the contraction
$i_{v}(*\phi)$ restricts to zero on $TP_{x}$. This leads to a \lq\lq Floer
type'' description of associative submanifolds. Let $P_{0}\subset M$ be some
compact submanifold
    and let ${\cal P}_{0}$ be a space of submanifolds close to $P_{0}$ in
a suitable sense. Then we can define a functional ${\cal F}$ on ${\cal P}_{0}$
by
    \begin{equation}   {\cal F}(P)= \int_{W} *\phi \end{equation}
    where $W$ is a $4$-chain in a small neighbourhood of $P_{0}$ in $M$ with
$\partial W= P-P_{0}$. The facts that $*\phi$ is closed and that we are only
working  with small deformations of $P_{0}$ means that this functional is
well-defined, independent of the choice of $W$. The derivative of ${\cal
F}$ at $P$ is given by $$
    \int_{P} i_{v} *\phi , $$
    where $v$ is a variation vector field, and we see that the associative
submanifolds in ${\cal P}_{0}$ are exactly the critical points of ${\cal
F}$. Globally, on a whole space of submanifolds ${\cal P}$, we do not usually
get a well-defined functional but we have a well-defined closed $1$-form $d{\cal
F}$ on ${\cal P}$ whose zeros are associative submanifolds,  in the familar
way in Floer-type theories. The linearisation of the associative condition
is an elliptic differential operator of index zero acting on sections of
the normal bundle. (The fact that the index is zero is a consequence of the
variational description, which implies that the linearisation is a self-adjoint
operator.)
    
    The discussion above applies to any $G_{2}$ structure for which the $4$-form
$*\phi$  is closed. The condition that the 3-form $\phi$ is closed enters
in the calibrated theory: it means that a compact  associative submanifold
is absolutely volume minimising in its homology class with volume the homological
invariant $\langle [\phi],[P]\rangle$.  In this case we get at least some
partial compactness properties of the set of associative submanifolds in
a fixed homology class.
    A longstanding theme in the literature is the possibility of developing
an enumerative theory, \lq\lq counting'' associative submanifolds, or---more
ambitiously---defining Floer homology groups. The fundamental difficulty
in doing this comes from the possible formation of singularities and corresponding
failure of compactness. A  detailed understanding of this seems a long way
off but there is a standard conjecture in the field, due to Joyce \cite{kn:J},
which suggests that it may only be  necessary to consider three phenomena. (See also the discussion in \cite{kn:DW}.)
 We consider a generic $1$-parameter family $\times_{t}$, for $t\in [-1,1]$
 on a compact $7$-manifold $M$ and fix a $3$-dimensional homology class in $M$.
We leave imprecise the exact meaning of generic and the exact conditions
imposed on the structures. One expects that for all but a discrete set $S\subset
[-1,1]$ of parameter values $t$ there is  a finite set of associative submanifolds
in the given homology class  and that as we vary $t$ in $[-1,1]\setminus
S$ a signed count of these will be locally constant. The question is, what
singularities can develop at the exceptional parameter values $t\in S$? Joyce
conjectures that, in generic $1$-parameter families, one will only encounter
phenomena which we will refer to in this article as:
    \begin{enumerate}
    \item \lq\lq multiple covers'';
    \item \lq\lq crossing'';
    \item \lq\lq surgery triples''.
    \end{enumerate}
    
   In (1) we have in mind a situation where a $1$-parameter family of embedded
submanifolds  $\iota_{t}:\Pi\rightarrow M$ converge as $t\rightarrow 0$ to
a map $\iota_{0}:\Pi\rightarrow M$ which is a covering (possibly branched)
of its image.  Such behaviour has been extensively studied for pseudo-holomorphic
curves and there is some work of Doan and Walpuski in the case of associative submanifolds \cite{kn:DW}, 
but the theory has not yet been developed very far so we will ignore
this multiple cover phenomenon in this article.
   
   \
   
   For (2), observe that for  dimensional reasons we expect that generically
associative submanifolds do not intersect (or self-intersect). However in
a generic $1$-parameter family we can expect to see associative submanifolds
$P_{t}, Q_{t}$ which intersect  at some parameter value $t=0$. It was predicted
by Joyce and confirmed by Nordstr\"om (in unpublished work, to appear) that in this situation (if the intersection in the family is transverse)
   there will be another family $(P\sharp Q)_{t}$ of associative submanifolds,
defined either for $t>0$ or for $t<0$ (but not both) described topologically by smoothing
the singular union $P_{0}\cup Q_{0}$ into a connected sum. The differential
geometric model near the intersection point is given by a \lq\lq Lawlor neck''.
   Thus is in this situation a straightforward count of associative submanifolds in the homology class of $(P\sharp Q)_{t}$ will change at $t=0$. 
   
   \
   
   In (3) the differential geometric model is provided by families of special
Lagrangian submanifolds in $\bC^{3}$ found by Harvey and Lawson \cite{kn:HL}.
 Take standard complex co-ordinates $z_{1}, z_{2}, z_{3}$ and for $s\geq
0$ define
   $$L_{1}^{s}= \{ (z_{1}, z_{2}, z_{3}): {\rm Im}(z_{1} z_{2} z_{3})=0,
{\rm Re}(z_{1}z_{2}z_{3})\geq 0, \vert z_{1}\vert^{2}-\vert z_{2}\vert^{2}=s;
\vert z_{2}\vert^{2}=\vert z_{3}\vert^{2} \}.
   $$
 When $s=0$ this is the cone over the standard torus $T^{2}\subset S^{5}$.
For $s>0$ we get a special Lagrangian submanifold in $\bC^{3}$ which, for
a standard matching of the structures, 
  is an associative submanifold in $\bC^{3}\times \bR$. Topologically, this
submanifold is obtained from the cone by cutting out a neighbourhood of the
vertex and gluing $D^{2}\times S^{1}$ to the resulting $T^{2}$ boundary.
Permuting the co-ordinates gives similar families
  $L_{2}^{s}, L_{3}^{s}$, both equal to the cone when $s=0$. Topologically,
we obtain $L_{2}^{s}, L_{3}^{s}$ from $L_{1}^{s}$ (all for $s>0$) by performing
Dehn surgeries on the circle in $L_{1}^{s}$ formed by the core of $D^{2}\times
S^{1}$.  The manifolds form a \lq\lq surgery triple'', with a cyclic symmetry
between them, differing by the way in which $D^{2}\times S^{1}$ is attached
to $T^{2}$. 
  
  Returning to the compact $7$-manifold $M$ with a $1$-parameter family of
structures $\times_{t}$, Joyce explains in \cite{kn:J} that one could encounter
a singular associative $P$ at parameter value $t=0$ with a singularity modelled
on the cone as above. For $t<0$ there could be a an associative submanifold
$P_{1}^{t}$, locally modelled on $L_{1}^{-t}$, and for $t>0$ a pair of associative
submanifolds $P_{2}^{t}, P_{3}^{t}$, locally modelled on $L_{2}^{t}, L_{3}^{t}$.
Or it could happen that $P_{2}^{t}, P_{3}^{t}$ exist for $t<0$ and $P_{1}^{t}$ exists
for $t>0$.  But either way a straightforward count of associative submanifolds
will change across $t=0$.

  \section{Kovalev-Lefschetz fibrations}
  \subsection{Topology}
  We will use the non-standard terminology that the {\it K3 manifold} is
the oriented differentiable 4-manifold $X$ underlying any complex K3 surface.
We recall that, with its cup product form, $H^{2}(X;\bR)$ is isomorphic to
$\bR^{3,19}$. We write $\Lambda_{X}\subset \bR^{3,19}$ for the integer
lattice and $O(\Lambda_{X})\subset O(3,19)$ for its automorphism group.  We call a class $\alpha$ in $\Lambda_{X}$ with $\alpha^{2}=-2$ a
{\it $-2$-class}. We will frequently use Poincar\'e duality to identify $H_{2}(X)$
with $H^{2}(X)$.
  \newtheorem{defn}{Definition}
  \newcommand{\tL}{\tilde{L}}
  \newcommand{\bZ}{{\bf Z}}
  \begin{defn}
  A topological Kovalev-Lefschetz (KL) fibration consists of data $(M,\Phi, N,L,\tL,\pi)$
where:
  \begin{itemize} \item $M$ is a compact $7$-manifold, and $\Phi$ is a class
in $H^{3}(M,\bR)$;
  \item $N$ is a compact oriented $3$-manifold, $L\subset N$ is a link and
$\pi:M\rightarrow N$ is a smooth map;
  \item $\tL\subset M$ is a submanifold and $\pi$ restricts to a diffeomorphism
from $\tL$ to $L$;
  \item at each point of $M\setminus \tL$ the derivative of $\pi$ is surjective;
  \item around a point of $\tL$ and the corresponding point of $L$ there
are co-ordinates in which $\pi$ is given by the model $\pi_{0}:\bC^{3}\times
\bR\rightarrow \bC\times \bR$
  \begin{equation}   \pi_{0} (z_{1}, z_{2}, z_{3}, t)= (z_{1}^{2}+z_{2}^{2}+z_{3}^{2};
t) \end{equation}
  \item each fibre of $\pi$ over points of $N\setminus L$ is diffeomorphic
to the K3 manifold $X$.
  \end{itemize}
  \end{defn}
  
  It is sometimes convenient to regard $N$ as an orbifold, with local orbifold degree $2$ at points of $L$.

  \

  Given this structure we obtain a flat vector bundle $E_{0}$  over $N\setminus L$ with fibre $\bR^{3,19}$ and structure group $O(\Lambda_{X})$ given by the cohomology along the fibres.
The monodromy of this bundle around a small loop about $L$ is of order $2$,
defined by reflection in a \lq\lq vanishing cycle''---which is a $-2$ class
in $H^{2}(X)$. We can extend this flat vector  bundle to a flat orbifold
vector bundle $E$ over $N$ and we have a corresponding sheaf $\underline{E}$
over $N$ of locally constant sections. The Leray spectral sequence gives
an exact sequence
  \begin{equation}   0\rightarrow  \bR=H^{3}(N;\bR)\rightarrow  H^{3}(M;\bR)\rightarrow
H^{1}(N;\underline{E})\rightarrow 0. \end{equation}
  The class $\Phi$ defines a lift of $E$ to an affine orbifold bundle $E^{+}$
over $N$. This is equivalent to saying that we have  flat  bundle $E_{0}^{+}
$ over $N\setminus L$ with structure group
  $A$, where $A$ is the affine extension
  $$    (\bR^{3,19},+)\rightarrow A \rightarrow O(\Lambda_{X}), $$
  and the monodromy around each component of $L$ maps to a reflection in
$O(\Lambda_{X})$. 
  To summarise, from any topological KL fibration as above
we can obtain data
  $(N,L,E^{+})$. The theme of this article is that we can study such data
in 3-dimensions, independent of the 7-dimensional picture which motivated
it.
  
  \

  In this article the dual of the sequence (3) will be important. Recall the
standard notion in algebraic topology of homology in a local co-efficient
system. In our case we have a local co-efficient system  $E_{0}$ over $N\setminus
L$. The group $H_{1}(N\setminus L, E_{0})$ can be defined by 1-cycles as
follows. We consider an oriented graph $\Gamma$ embedded in $N\setminus L$.
Each edge $\gamma$ of the graph is labelled by a constant section $\alpha_{\gamma}$
of $E_{0}$ over $\gamma$ and reversing the orientation of $\gamma$ is equivalent
to changing the sign of the label. At a vertex of the graph we have
  \begin{equation} \sum_{\gamma} \pm \alpha_{\gamma}=0 \end{equation}
  where the sum runs over edges $\gamma$ incident to the vertex and the sign
$\pm$ is determined by ingoing/outgoing orientation.  To obtain the homology
group $H_{1}(N\setminus L,E_{0})$ we divide the group of $1$-cycles by a
subgroup of boundaries, defined in a similar fashion.  
  
  To bring in the link $L$ we make one change to the above recipe. We take
graphs in $N$ and we allow an edge to terminate on a point of $L$ provided
that its label is a multiple of the corresponding vanishing cycle. (To be  more precise, the vanishing cycle is defined up to sign in the fibres of $E_{0}$ near $L$ and we allow the label to take either sign.)
In this
fashion we obtain a group which we denote by $H_{1}(N,E)$ which fits into
an
   exact sequence \begin{equation} 0\rightarrow H_{1}(N;E)\rightarrow H_{3}(M;\bR)\rightarrow
H_{3}(M;\bR)=\bR\rightarrow 0 \end{equation} dual to (3). 
  
  The class $\Phi\in H^{3}(M)$ gives a linear map from $H_{3}(M)$ to $ \bR$
and hence a linear map $\chi: H_{1}(N;E)\rightarrow \bR$. This can be described
as follows. Choose any smooth orbifold section $u$ of the affine bundle $E^{+}$.
 That is, in an orbifold chart around a point of $L$ the section is given
by a $\bZ/2$-equivariant map to $\bR^{3,19}$, where $\bZ/2$ acts on $\bR^{3,19}$
by the local monodromy.  Let $\Gamma$ be a labelled graph as above. For each
edge $\gamma$ of $\Gamma$ with label $\alpha_{\gamma}$ and running from $p$
to $q$ we can define
  \begin{equation}  \langle \gamma, u\rangle=  \alpha_{\gamma}(u(p)- u(q)).
\end{equation}
  Then the map $\chi: H_{1}(N;E)\rightarrow \bR$ defined by $\Phi$ is induced at
the chain level by
$ \Gamma\mapsto \sum_{\gamma} \langle \gamma, u\rangle$.

  \subsection{Adiabatic $G_{2}$ structures}
  
  We recall some material from \cite{kn:D1}.  One standard model for the
cross-product on $\bR^{7}$ is obtained by writing $\bR^{7}= \bR^{4}\oplus
\bR^{3}$ where $\bR^{4}$ is taken with its standard orientation and Euclidean
structure and $\bR^{3}$ is identified the $3$-dimensional space of self-dual
$2$-forms on $\bR^{4}$. It is more convenient to work with the $3$-form (which
determines the cross-product). If $\omega_{1}, \omega_{2},\omega_{3}$ is
a standard basis for $\Lambda^{+}$  and $y_{i}$ are co-ordinates on $\bR^{3}$,
the model $3$-form on $\bR^{7}$ is
  $ \phi_{0}= -\sum \omega_{i} dy_{i}+ dy_{1}dy_{2}dy_{3}$.
  (Here we are using a different sign convention from \cite{kn:D1}, which fits better
with the literature.)
  
  We define  a {\it hyperk\"ahler structure} on the K3 manifold X to be a
 triple of closed $2$-forms $\Omega_{1},\Omega_{2}, \Omega_{3}$ such that

  $$  \Omega_{i}\wedge \Omega_{j}= a_{ij} {\rm vol}_{X}, $$
  for some volume form ${\rm vol}_{X}$ on $X$ and a constant positive definite
matrix $(a_{ij})$. It is more standard to require that $a_{ij}=\delta_{ij}$
but this can always be achieved by a change of basis and the extra freedom
will be convenient.   
   If we have such a structure then we get a torsion-free $G_{2}$-structure
on $X\times \bR^{3}$  with $3$-form \begin{equation} -\sum \Omega_{i} dy_{i}+ dy_{1}dy_{2}dy_{3}. \end{equation}
The deep fact we need is the {\it Torelli theorem} for K3 surfaces. In our
set-up this can be stated as follows. Suppose that $h_{1}, h_{2}, h_{3}\in
\bR^{3,19}=H^{2}(X;\bR)$ span a maximal subspace $H$ in $\bR^{3,19}$. Suppose
in addition that there is no class $\alpha$ in $\Lambda_{X}$ with $\alpha^{2}=-2$
which is orthogonal to $H$. Then there is a hyperk\"ahler structure $(\Omega_{i})$
on $X$ with $[\Omega_{i}]=h_{i}$ and this is unique up to the action of diffeomorphisms
of $X$ which act trivially on $H^{2}(X)$.
   
   We now return to our topological KL fibration and affine
bundle $E^{+}$ over the $3$-manifold $N$ and  let $u$ be a section of $E^{+}$.
We consider first the situation over a co-ordinate neighbourhood $B\subset
N\setminus L$ with co-ordinates $y_{i}$. Over $B$ the section $u$ is given
by a map  $u_{B}:B \rightarrow \bR^{3,19}$. We say that $u_{B}$ is positive if at each point $b$
of $B$ the image of the derivative of $u_{B}$ is a maximal positive subspace
$H_{b}$ in $\bR^{3,19}$. We also assume that 
   there is no $-2$ class orthogonal to $H_{b}$.  Then we can apply the Torelli
theorem with $h_{i}$ the image of $\partial_{y_{i}}$ under the derivative
of $u_{B}$.  So at each point $b$ of $B$ we get a hyperk\"ahler structure $(\Omega_{i})$
on the fibre $\pi^{-1}(b)\subset M$. A procedure for   choosing forms $\tilde{\Omega}_{i}$
on $\pi^{-1}(B)$ which restrict to $\Omega_{i}$ on the fibres is explained
in \cite{kn:D1}.  We say that $u_{B}$ is  {\it
maximal positive } if it is positive and the image is a \lq\lq maximal submanifold'' of $\bR^{3,19}$---i.e. the image satisfies the Euler-Lagrange equation
associated to the volume functional, just as for minimal submanifolds in Euclidean spaces.  Write the induced volume form on $B$ as
$\lambda dy_{1}dy_{2}dy_{3}$ and introduce a positive parameter $\epsilon$.
Then we call the $3$-form on $\pi^{-1}(B)\subset M$
   \begin{equation} \phi_{\epsilon}= -\epsilon \sum \tilde{\Omega}_{i} dy_{i}+ \lambda
dy_{1}dy_{2}dy_{3}\end{equation}
   an {\it adiabatic solution} to the torsion-free $G_{2}$-equations. (This is independent of the choice of local co-ordinates.)  For
 a fixed $\epsilon$ this is not an exact solution but, roughly speaking,
when $\epsilon$ is small it can be deformed slightly to an exact solution $\tilde{\phi}_{\epsilon}$. (For more
precise statements we refer to \cite{kn:D1}.) In the metric induced by $\phi_{\epsilon}$ the diameter of the fibre is $O(\epsilon^{1/2})$ while the diameter of $\pi^{-1}(B)$ is $O(1)$, so as $\epsilon\rightarrow 0$ we are studying a collapsing family of metrics (in the same spirit as the Calabi-Yau metrics discussed by Tosatti in this volume \cite{kn:T}).
   
   The most technical issue here is the behaviour required of the section $u$ near a point $p$ of the link $L$, which we now describe.
\begin{itemize}
\item There is an orthogonal decomposition \begin{equation} \bR^{3,19}= H\oplus \bR \delta \oplus K \end{equation}
where $H$ is a maximal positive subspace and $\delta$ is the vanishing cycle. We require that $\pm \delta$ are the only $-2$ classes in $H\oplus \bR \delta$.
\item There  are local co-ordinates $(z,t)$ centred at $p$ (with $z\in \bC$ and $t\in \bR$) so that $L$ is locally defined by $z=0$ and the section $u$ is given by a multivalued map $(I(z,t), f(z,t), F(z,t))$ with respect to the decomposition of $\bR^{3,19}$ as follows.
\begin{enumerate} \item $I$ is the identity under a linear identification $H=\bC\oplus \bR$;
\item $F$ is a smooth single-valued function which vanishes along with its first derivatives at the origin;
\item $f$ is a 2-valued function (defined up to sign) with
\begin{equation} f(z,t) = {\rm Re}(b(t) z^{3/2}) + \eta(z,t) \end{equation}
where $b$ is smooth and nowhere zero and $\eta(z,t)=O(\vert z\vert ^{5/2})$. More precisely we require that all $t$ derivatives of $\eta$ are $O(\vert z\vert^{5/2})$ and the first derivatives of $\eta$ in the $z$ factor are $O(\vert z\vert^{3/2})$. 
\end{enumerate}
\end{itemize}

In a standard way, we can avoid the use of multivalued functions by passing
to an orbifold chart, as explained in \cite{kn:D1}. We define a {\it branched maximal positive section} of $E^{+}$ to be  a section which satisfies these conditions near $L$ and the conditions described before in local coordinates $B\subset N\setminus L$.

\

The conclusion is that we have a notion of an {\it adiabatic $G_{2}$-structure} which is a quadruple $(N,L, E_{+}, u)$ where $E_{+}$ is a flat affine orbifold bundle over $N$ with monodromy around $L$ given by reflections and $u$ is a branched maximal positive section.

\

We should mention that, at the time of writing,  there are no real examples of these structures known and it has not been shown that they lead to  collapsing families of torsion-free $G_{2}$-structures $\tilde{\phi}_{\epsilon}$ on  compact $7$-manifolds.  But one can hope that this will change. 

\section{Adiabatic associative submanifolds}

Let $(\Omega_{i})$ be a hyperk\"ahler structure on $X$. For each nonzero vector $v\in \bR^{3}$ there is a complex structure $J_{v}$ on $X$ for which $\sum v_{i}\Omega_{i}$ is a K\"ahler form.  If $\Sigma\subset X$ is a smooth complex curve with respect to the complex structure $J_{v}$ then it is easy to check that $\Sigma \times \bR v$ is an associative submanifold in the product $X\times \bR^{3}$ with the $G_{2}$-structure (7). Let $\alpha$ be a class in the integer lattice $\Lambda_{X}\subset \bR^{3,19}$ with $\alpha^{2}\geq -2$ and let $p(\alpha)\in H$ be the orthogonal projection to the subspace $H$ spanned by $[\Omega_{i}]$ so $p(\alpha)=\sum v_{i}[\Omega_{i}]$ for a vector $v$ in $\bR^{3}$.
The vector $v$ is not zero: if $\alpha^{2}\geq 0$ this follows from the fact that $H$ is a maximal positive subspace and if $\alpha^{2}=-2$ it is the condition arising in the Torelli theorem. By construction, $\alpha$ has type $(1,1)$ with respect to the complex structure $J_{v}$ and it corresponds   to a holomorphic line bundle $L$ of positive degree. The  Riemann-Roch theorem implies that
${\rm dim}\  H^{0}(L)\geq \alpha^{2}/2 +2\geq 1$,  so $\alpha$ is represented by a $J_{v}$-complex curve. In this article we restrict attention to the case when $\alpha^{2}=-2$, that is, a $-2$ class. Then standard theory gives that there is a {\it unique} $J_{v}$ complex curve $\Sigma$  in the class $\alpha$ and this is \lq\lq generically'' a smooth embedded $2$-sphere. More precisely, we will say that the pair $(\alpha,H)$ is {\it irreducible} if we cannot write
$\alpha$ as a sum of $-2$ classes $\alpha_{i}$ with $p(\alpha_{i})=\lambda_{i} p(\alpha)$ for $\lambda_{i}>0$. Then if $(\alpha,H)$ is irreducible the $J_{v}$ complex representative is a smooth $2$-sphere. If  $(\alpha,H)$ is not irreducible  the representative could have a number of components,  as we will discuss further in 6.4 below.

Let $\phi_{\epsilon}$ be the adiabatic 3-form in (8) over $\pi^{-1}(B)\subset M$. To simplify the formulae we work in local co-ordinates centred at a point $b$ in $B$, chosen so that $a_{ij}=\delta_{ij}$ at $b$ and the integral of the volume form ${\rm vol}_{X}$ is $1$. So, at the point $b$,
$$  \int_{X} \Omega_{i}\wedge \Omega_{j}= \delta_{ij} $$
and $\lambda=1$. Then at this point 
\begin{equation} *\phi_{\epsilon}= \epsilon  \sum \tilde{\Omega}_{i} dy_{j}dy_{k} + \epsilon^{2} {\rm vol}_{X} , \end{equation} where the sum runs over cyclic permutations of $(123)$. 
Thus the family of $4$-forms $\epsilon^{-1} *\phi_{\epsilon}$ has a well-defined limit $\Psi= \sum \Omega_{i}dy_{j}dy_{k}$ as $\epsilon\rightarrow 0$. Of course $\Psi$ is defined over all of $B$, with a more complicated formula in a general co-ordinate system.
Recall that the associative condition for a submanifold in the $G_{2}$-structure $\phi_{\epsilon}$ is that $i_{v}*\phi_{\epsilon}$ restricts to $0$ for all tangent vectors $v$ in the $7$-manifold. Thus, while the $4$-form $\Psi$ does  {\it not} correspond to a $G_{2}$-structure we can consider the same condition. We say that a submanifold $P$ is \lq\lq $\Psi$-associative'' if  $i_{v}\Psi$ restricts to zero on $P$ for all $v\in TM$.

Let
$\alpha$ be a $-2$ class. In a  trivialisation of the bundle $E^{+}$ over $B$ the section $u$ is given by a map $u_{B}:B\rightarrow \bR^{3,19}$ and we have a function $h_{\alpha}$ on $B$ defined by
$$   h_{\alpha}(y)= \langle \alpha, u_{B}(y)\rangle . $$
We also have a Riemannian metric on $B$ induced by the embedding in $\bR^{3,19}$. Let   $\gamma:(a,b)\rightarrow B$ be a {\it gradient flowline}  for the function $h_{\alpha}$, {\it i.e.} 
$$    \gamma'(s)= \left( {\rm grad} h_{\alpha}\right)_{\gamma(s)}. $$
Assume that at each point $\gamma(s)$ the pair $(\alpha,H_{\gamma(s)})$ is irreducible,  where $H_{\gamma(s)}$ is the positive subspace in $\bR^{3,19}$ defined by the image of $du_{B}$. We can apply the preceding discussion at each point  $\gamma(s)$ and we see that there is a unique  embedded $2$-sphere  $\Sigma_{s}\subset \pi^{-1}(\gamma(s))$ in the class $\alpha$ which is complex with respect  to the complex structure defined by the velocity vector $\gamma'(s)$. Thus we obtain  a $3$-dimensional submanifold $P_{\gamma}\subset \pi^{-1}(B)$ which fibres over the flowline $\gamma$ with fibre over $\gamma(s)$ the $2$-sphere $\Sigma_{s}$. Then we have
\newtheorem{prop}{Proposition}
\begin{prop}
The submanifold $P_{\gamma}$ is   $\Psi$-associative.
\end{prop}
To see this we can work in adapted co-ordinates at a point in $B$ as above and assume that the gradient of $h_{\alpha}$ is a  multiple of $\partial_{y_{1}}$. From the nature of $P_{\gamma}$ the condition that $i_{v}\Psi$ restricts to $0$ is automatic if $v$ is a tangent vector along the fibres or a lift of $\partial_{y_{1}}$.  Taking  lifts of $\partial_{y_{2}}, \partial_{y_{3}}$ the $\Psi$-associative condition is just that $\Omega_{2}$ and $\Omega_{3}$ vanish on the fibre of $P_{\gamma}$,  which is just the condition that this fibre is a complex curve with respect to the complex structure corresponding to $\partial_{y_{1}}$.

\section{Calibrated 1-cycles}

We will now formulate the main definition of this article. Let $(N,L,E^{+}, u)$ be an adiabatic $G_{2}$ structure. We want to consider two situations: \lq\lq orbits'' and \lq\lq graphs''.

\

{\bf Orbits}

Let $\Gamma\subset N\setminus L$ be an oriented embedded circle.
Suppose that there is a constant $-2$ section $\alpha$ of the bundle $E_{0}$ over $\Gamma$, so the pair $(\Gamma,\alpha)$ is a cycle for the group $H_{1}(N,E)$ discussed in 3.1. Locally, we have a function $h_{\alpha}$ as considered above. This function is not globally defined but its gradient vector field is defined on a neighbourhood of $\Gamma$ in $N\setminus L$. We say that $(\Gamma,\alpha)$ is {\it a gradient orbit} if it is an integral curve of this vector field (compatible with the orientation of $\Gamma$ in the obvious way).  We also require that for each point $q$ on $\Gamma$ the pair $(\alpha, H_{q})$ is irreducible,  where $H_{q}$ is the positive subspace defined by the derivative of $u$ at $q$. 

\

{\bf Graphs}

Let $\hat{\Gamma}$ be an oriented graph with each vertex having valence $1$ or $3$.  Let $\iota:\hat{\Gamma}\rightarrow N$ be a continuous embedding which is smooth on each edge and which maps the vertices of valence $1$ to $L$ and all other points of $\hat{\Gamma}$ to $N\setminus L$. We call the image $\Gamma$ of $\iota$ an embedding of $\hat{\Gamma}$. (Thus $\Gamma$ does not depend on the  parametrisation $\iota$.) Suppose that for each edge $\gamma$ of $\Gamma$ there is a label by a constant $-2$ section $\alpha_{\gamma}$ of $E$,  that the label of the edge containing a vertex of valence $1$ is the corresponding vanishing cycle and that the condition (4) holds at each vertex of valence $3$. So $(\Gamma, \{\alpha_{\gamma}\})$  is a cycle for the group $H_{1}(N,E)$.
We say that  $(\Gamma, \{\alpha_{\gamma}\})$ is a {\it gradient graph} if each edge $\gamma$ is an integral curve of the gradient vector field defined by $\alpha_{\gamma}$ (compatible with the orientations). We also require that the  irreducible condition holds at each point of $\Gamma$. (At a vertex $q$ of valence $3$ this means that $(\alpha_{i}, H_{q})$ is irreducible for each of the $-2$ classes $\alpha_{1}, \alpha_{2}, \alpha_{3}$ labelling the three edges.) Note that the sign ambiguity in the vanishing cycle near to $L$ is taken care of by the orientation of the graph.

\

\begin{figure}[t]
\centering
\includegraphics[scale=.75]{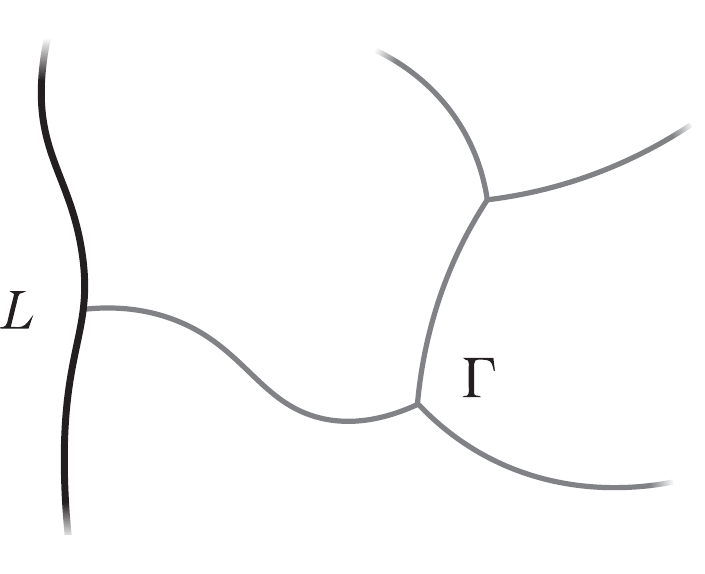}
\caption{Local diagram of gradient cycle $\Gamma$ which has one univalent vertex terminating at the link $L$.}
\label{fig:calibrated-cycle}
\end{figure}

In either case we call $\Gamma$ (with the appropriate labels) a  {\it gradient cycle}. See Figure \ref{fig:calibrated-cycle} for an example. 

\

{\bf Discussion}
\begin{enumerate}
\item 

We could vary the definition of gradient cycles in a  number of ways. We could
\begin{enumerate}
\item allow vertices of higher valence;
\item relax the condition that $\iota$ is an embedding; 
\item relax the irreducible requirement  on $(\alpha, H_{q})$;
\item allow $\Gamma$ to meet $L$ in interior points of edges.
\end{enumerate}
But the definition we have given simplifies some statements and one expects that for generic sections $u$  the gradient cycles obtained under such  weaker hypotheses will satisfy our stronger conditions. The topic of Section 6 below is to study situations where these various conditions fail, in $1$-parameter families.

\item

The fundamental idea underlying these definitions is that if the adiabatic $G_{2}$-structure $(N,L,E^{+},u)$ corresponds to a $1$-parameter family $\tilde{\phi}_{\epsilon}$ of KL fibred torsion-free $G_{2}$-structures on a $7$-manifold $M$ then for small $\epsilon$ a gradient $1$-cycle $\Gamma$ should yield an associative submanifold $\tilde{P}(\Gamma)$ in $(M,\tilde{\phi}_{\epsilon})$. In the case of orbits it is relatively straightforward to prove such a result. First, at the topological level, since the irreducible condition holds over $\Gamma$ we clearly have a compact $3$-manifold $P(\Gamma)\subset M$ fibered over $\Gamma$ with 2-sphere fibres (and in fact $P(\Gamma)$ is diffeomorphic to $S^{2}\times S^{1}$). Then the analysis problem is to show that this can be deformed into an associative submanifold.

In the case of a graph $\Gamma$ the problem is harder but, as a first step, we will construct a topological model:  a submanifold $P(\Gamma)\subset M$.  For a vertex of valence $1$ on an edge $\gamma$ which terminates  at a point $p$ in $L$ the construction is essentially the well-known \lq\lq thimble'' of a Lefschetz fibration. For the topological discussion in our model around $p$ we can deform $\gamma$ to be the positive real axis in $\bC\subset \bC\times \bR$. Then the thimble is simply given by $\bR^{3}\subset \bC^{3}\subset \bC^{3}\times \bR$ and the fibre of $\pi_{0}$ over a point $\eta>0$  is the $2$-sphere $\{(z_{1}, z_{2}, z_{3}): z_{i}\in \bR, \sum z_{i}^{2}=\eta\}$.

The case of a vertex of valence $3$ is more interesting. We can suppose that all orientations are outgoing so we have three $-2$ classes $\alpha_{1}, \alpha_{2}, \alpha_{3}$ with $\alpha_{1}+\alpha_{2}+\alpha_{3}=0$. This implies that $\alpha_{i}.\alpha_{j}=1$. In the hyperk\"ahler structure on the fibre $X$ these classes are represented by embedded spheres $\Sigma_{1}, \Sigma_{2}, \Sigma_{3}$ which are complex with respect to three complex structures $J_{1}, J_{2}, J_{3}$. Standard theory of K3 surfaces shows that this configuration of embedded spheres is unique up to diffeomorphisms of $X$, so for the topological discussion we can take any convenient model. The model we use involves the non-compact manifold $X_{0}$ obtained as either a smoothing or resolution of the $A_{2}$-singularity which, again by standard theory, can be emebedded in the K3 manifold $X$. We use  the \lq\lq Gibbons-Hawking'' description of $X_{0}$ (but for our present  purposes only at the topological level).  For this we take three points $A,B,C$ in $\bR^{3}$ and construct an $S^{1}$-bundle $Z\rightarrow \bR^{3}\setminus \{A,B,C\}$ with Chern class $1$ on small spheres around the three points. The restriction of $Z$ to such a sphere gives the Hopf fibration of $S^{3}$ so we can complete $Z$ to a $4$-manifold $X_{0}$ which has  a circle action with three fixed points and the quotient by the action is a smooth map $\mu:X_{0}\rightarrow \bR^{3}$. The pre-image by $\mu$ of the line segment from $A$ to $B$ is a $2$-sphere in $X_{0}$ with self-intersection $-2$ and taking the three sides of the triangle $ABC$ we get a configuration of spheres of the desired kind.  We can suppose that $A,B$ and $C$ lie in the standard plane $\bR^{2}\subset \bR^{3}$. Let $v_{1}=A-B, v_{2}=B-C, v_{3}=C-A$,  so $v_{i}$ are vectors in $\bR^{2}$ with  $\sum v_{i}=0$. Let $Y\subset \bR^{2}$ be the union of the three rays $\bR^{+} v_{i}$ and let $\Omega\subset \bR^{2} $ be a  thickening of $Y$, as in Figure \ref{fig:y-graph}, with three boundary components, asymptotic  to the three pairs of rays. Let $F:\Omega\rightarrow \bR^{2}$ be a smooth map which takes the three boundary components to the points $A,B,C$ with the obvious ordering (so, for example,  the boundary component asymptotic to $v_{1}$ and $v_{3}$ is mapped to $A$). We  choose the map $F$
so that outside a large ball it maps into the three edges of the triangle, again in the obvious way. The graph of $F$ is a surface $S$ with boundary in $\bR^{2}\times \bR^{2}\subset \bR^{3}\times \bR^{2}$ and the boundary of $S$ lies in $\{A,B,C\}\times \bR^{2}$. It follows that the preimage
$$P_{0}= (\mu\times {\rm id})^{-1}(S)\subset X_{0}\times \bR^{2}$$ is a $3$-dimensional submanifold of $X_{0}\times \bR^{2}$. By construction, $P_{0}$ has three ends, which are small deformations  of  $\Sigma_{i}\times \bR^{+}v_{i}$. Finally, taking $\bR^{2}\subset \bR^{3}$ we can regard $P_{0}$ as a  submanifold of $X_{0}\times \bR^{3}$.

\begin{figure}[t]
\centering
\includegraphics[scale=0.9]{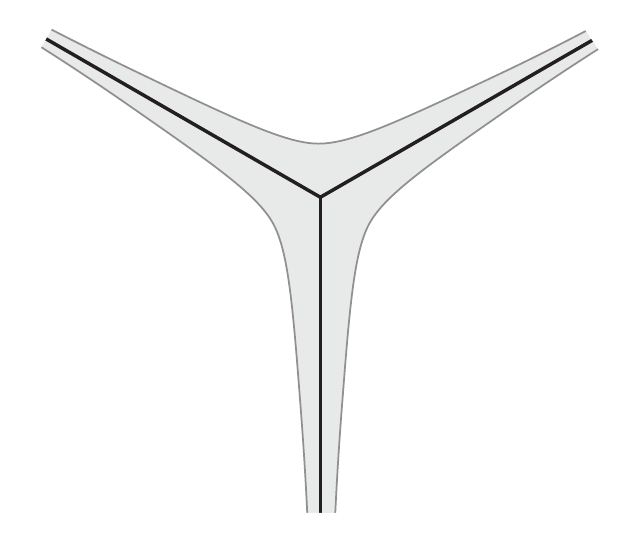}
\caption{Neighborhood of $Y\subset \bR^2$}
\label{fig:y-graph}
\end{figure}

There are many other ways of describing this submanifold $P_{0}\subset X_{0}$. The approach above has the advantage that it makes evident the symmetry between the three 2-spheres. 
Using this model at each vertex of valence $3$ we can construct a $3$-dimensional submanifold $P(\Gamma)\subset M$ for any graph $\Gamma$.

The problem of  constructing  a nearby associative submanifold $\tilde{P}(\Gamma)$ motivates the following conjecture.
\newtheorem{conj}{Conjecture}

\begin{conj} 
Let $\alpha_{1}, \alpha_{2}, \alpha_{3}$ be $-2$ classes on the K3 manifold $X$ with $\alpha_{1}+\alpha_{2}+\alpha_{3}=0$. Let $\bR^{3}=H\subset H^{2}(X)$ be a maximal positive subspace corresponding to a hyperk\"ahler structure and $v_{i}$ be the projection of $\alpha_{i}$ to $H$. Assume that the $(\alpha_{i},H)$ are irreducible. Then there is an associative submanifold $\Pi\subset X\times \bR^{3}$ with three ends asymptotic to $\Sigma_{i} \times \bR^{+}v_{i}$ where $\Sigma_{i}$ is the complex curve representing $\alpha_{i}$,  for the complex structure defined by $v_{i}$, and $\Pi$ is unique up to the translations of $\bR^{3}$.
\end{conj} 

Since the vectors $v_{i}$ lie in a plane this associative submanifold should in fact  be a {\it special Lagrangian} submanifold in $X\times \bC$.

\item 

Gradient cycles are the zeros of  a closed $1$-form on a suitable infinite-dimensional space of cycles. This is simplest to set up in the case of orbits. Let $\Gamma_{0}$ be an embedded circle in $N\setminus L$ with constant section $\alpha$ over $\Gamma_{0}$ and let $V_{\alpha}$ be the corresponding gradient vector field, defined in a neighbourhood of $\Gamma_{0}$. Let $\sigma$ be the $2$-form on this neighbourhood given by the contraction of the volume form with $V_{\alpha}$. Locally we have a function $h_{\alpha}$ and $\sigma= * dh_{\alpha}$. Recall that $h_{\alpha}$ is obtained as the restriction of a linear function on $\bR^{3,19}$ to the image of $U$, which is a maximal submanifold. The maximal condition implies that the restriction of linear functions to the submanifod are harmonic, with respect to the induced metric. So $h_{\alpha}$ is a harmonic function and $\sigma$ is a closed $2$-form. For any $1$-cycle $\Gamma$  close to $\Gamma_{0}$ choose a $2$-chain $W$ of the obvious kind with $\partial W=\Gamma-\Gamma_{0}$ and define
$$  {\cal F}(\Gamma)= \int_{W} \sigma. $$
Then the gradient orbits in this neighbourhood are the critical points of this functional. The derivative is a well-defined $1$-form given by
$$    \delta {\cal F}= \int_{\Gamma}  i_{\xi} \sigma, $$
where $\xi$ is a variation vector field along $\Gamma$.

We can proceed in a similar way for graphs. For each edge $\gamma$ we consider a variation vector field $\xi_{\gamma}$. These satisfy a matching condition $  \xi_{\gamma_{i}}= \xi_{\gamma_{j}}$ at a vertex of valence $3$ and are tangent to $L$ at vertices of valence $1$.
We define a $1$-form by
$$  \delta{\cal F}= \sum_{\gamma} \int_{\gamma} i_{\xi_{\gamma}} \sigma_{\gamma}. $$ where $\sigma_{\gamma}$ is the closed 2-form in a neighbourhood of $\gamma$ defined as above.  The reader can check that this is a closed $1$-form on an infinite dimensional space of labelled embedded graphs, with zeros the gradient graphs. 

One can go on to develop at least some elements of a \lq\lq Floer theory'',  with chain complex generated  by gradient graphs and boundary map defined by suitable $2$-cycles in $N\times \bR$ (which should correspond to Cayley submanifolds in $M\times \bR$). In the same vein, one can develop  theory for adiabatic co-associative submanifolds in $M$,  but we will not go into these variants further in this article.

\

There is also an adiabatic analogue of the calibrated property. For a segment  $\gamma(s)$ with label $\alpha$  we define the weighted length to be
$$  \int \vert \gamma'(s) \vert \vert \nabla h_{\alpha} \vert_{\gamma(s)} ds. $$
Then  the weighted length of any cycle $\Gamma$ is bounded below by $\vert \chi [\Gamma]\vert $ where $\chi:H_{1}(E)\rightarrow \bR$ is the map  discussed in 3.1 above and equality holds if and only if $\pm \Gamma$ is a gradient cycle.  
\end{enumerate}

\newtheorem{lem}{Lemma}

\subsection{Related literature}
There is a web of connections between the  ideas that we discuss in this article and existing literature. To give a proper account of this would go far beyond the space available here and the authors' knowledge,  so we will just  indicate some of these connections.

\begin{itemize}
\item In  classical algebraic topology, graphs of gradient lines on a manifold can be used to describe cup products and higher operations in the framework of Morse Theory \cite{kn:CN}.
\item In symplectic topology, similar graphs appear in describing pseudo-holomorphic curves \cite{kn:FO}, \cite{kn:A} in various contexts. These include  relations to tropical geometry,  Calabi-Yau manifolds with Lagrangian torus fibrations, and the Strominger-Yau-Zaslow approach to mirror symmetry \cite{kn:GS}.
\item Invariants of $3$-manifolds, related to Chern-Simons theory,  obtained by counting graphs of gradient lines have been studied by Fukaya \cite{kn:F} and Watanabe \cite{kn:W}. This leads to connections with the theory of finite type invariants of 3-manifolds, see e.g. \cite{kn:L}.
\item In the more specific context of $G_{2}$-manifolds fibred over a $3$-dimensional base with K3 fibres, the basic idea in this article was discussed in \cite{kn:D1} for the particular case of an arc joining two components of the link. Some closely related ideas appear in the earlier paper \cite{kn:PW}, in the case where the fibre is a hyperk\"ahler  ALE 4-manifold. Similar ideas have been considered in the case of Special Lagrangian submanifolds in Calabi-Yau 3-folds fibred over a 2-dimensional base (\cite[Section 3]{kn:K}, \cite[Section 1.4]{kn:S}) and there are connections to the theory of Spectral Networks. A recent article \cite{kn:KL}  gives some explicit examples of fibrations of noncompact $G_{2}$-manifolds, and associative submanifolds.

\item In 3-manifold topology, Hutchings defined a \lq\lq periodic Floer homology'' using the integral curves of the gradient vector field of a circle-valued harmonic function, in the case when the vector field has no zeros \cite{kn:H}.   This was shown by Lee and Taubes to agree with Seiberg-Witten-Floer homology \cite{kn:LT}. In this case the $3$-manifold is necessarily a surface bundle over the circle. A general theory in  the case of a vector field with zeros does not seem to have been developed yet but it seems likely that the work of  Taubes \cite{kn:Ta} and Gerig \cite{kn:G} on holomorphic curves in near symplectic manifolds would be relevant to that.

\end{itemize}
 
\section{Theory of gradient cycles}

\subsection{Terminating manifolds}

\begin{figure}[t]
\centering
\includegraphics[scale=.35]{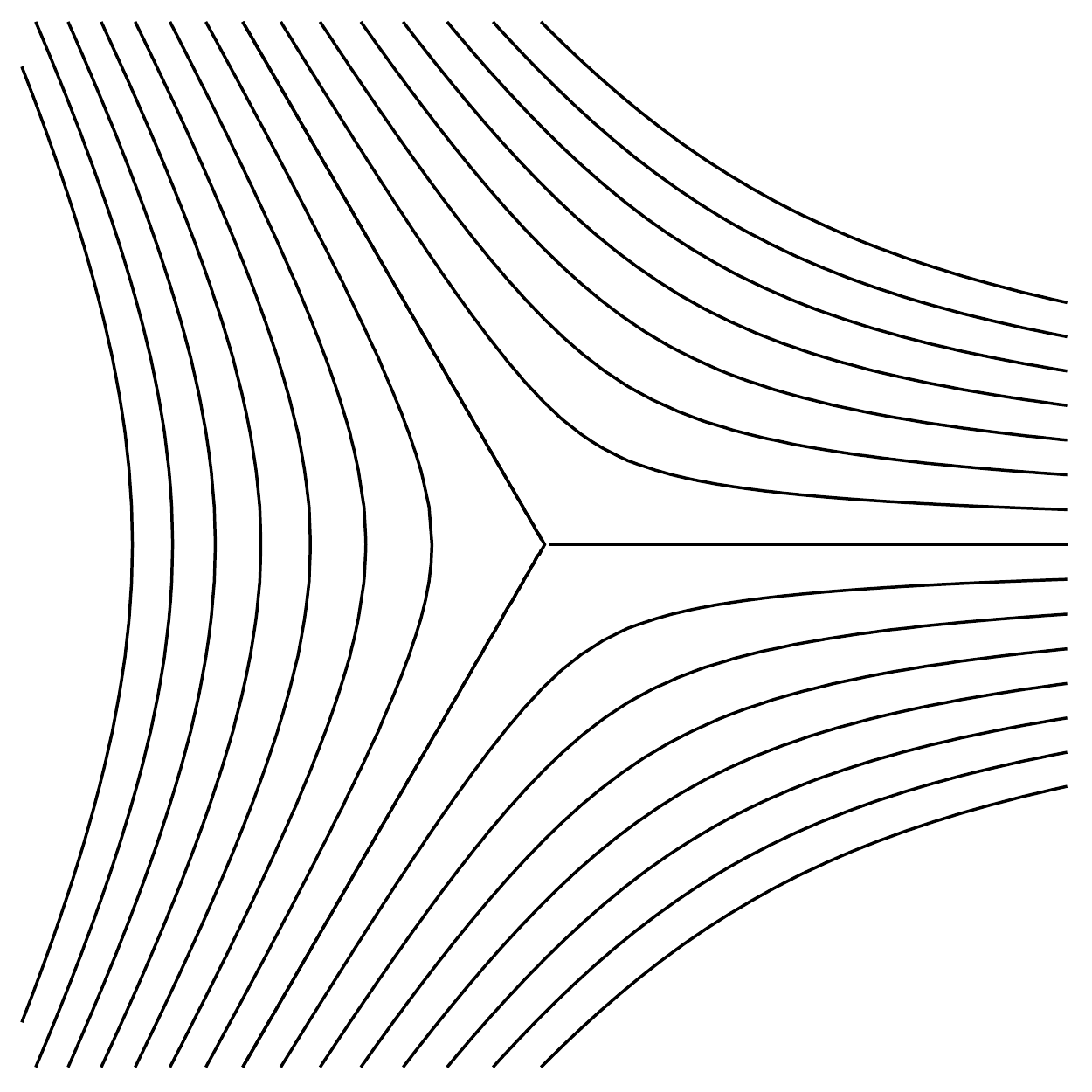}
\caption{Flow lines of the vector field $2 z^{1/2}$ on $\bC$.}
\label{fig:vector-field}
\end{figure}

In this subsection we  discuss  gradient flow lines corresponding to the vanishing cycle near the link $L$. We work in the unit ball $B^3$ in standard local co-ordinates $(z,t)$ as in 3.2, so we have maps $f:B^3\rightarrow \bR$ and $F:B^3\rightarrow \bR^{18}$ such that
$f(z,t)= {\rm Re}(b(t) z^{3/2}) +\eta(z,t)$, where $\eta, F$ satisfy the conditions stated in 3.2.  We have a (non-smooth) Riemannian metric $g$ on $B^3$ defined by the graph of $(f,F)$ in $\bR^{3,19}=\bC\times \bR \times \bR \times \bR^{18}$ and a vector field $V$ on $B^3$ given by the gradient of $f$ in the metric $g$.   Without loss of generality assume that $b(0)= 4/3$ and take the standard branch of $z^{3/2}$ near the positive real axis.  Let $V_{0}$ be the gradient vector field of ${\rm Re}( 4/3 z^{3/2})$ on $\bC\times \bR$ with respect to the {\it Euclidean} metric $g_{0}$. Thus $V_{0}= ( 2 z^{1/2}, 0)$ and there is an integral curve  $x_{0}(s)= ( s^{2}, 0)$ of $V_{0}$ that touches the origin. See Figure \ref{fig:vector-field}. Let $s_{0}$ be a small number to be fixed later.

\begin{prop}
There is a unique function $h:[0,s_{0})\rightarrow \bC\times \bR$ with $\vert h(s)\vert \leq C s^{3} $ and $\vert h'(s)\vert \leq C s^{2}$ such that $x_{0} +h $ is an integral curve of $V$.
\end{prop} 
 
To see this write $V=V_{0}+W$ so we want to solve the equation
\begin{equation}    h'= V_{0}(x_{0}+h)- V_{0}(x_{0}) + W(x_{0}+h) .\end{equation}
Straightforward calculations show that our hypotheses on $\eta, F$  imply that 
$$  \vert W(z,t)\vert = O(\vert z\vert^{3/2}+ \vert t\vert \vert z\vert^{1/2}), $$
with the corresponding estimates for the derivatives. We have an elementary estimate that if  $ \vert \zeta\vert\leq (1/2) \vert z\vert$ and if we set
$$  2(z+\zeta)^{1/2}-2z^{1/2}- z^{-1/2} \zeta = q(z,\zeta)$$
then $$\vert q(z,\zeta)\vert \leq C \vert \zeta\vert^{2} \vert z\vert^{-3/2}$$ and if $\vert \zeta_{1}\vert ,\vert \zeta_{2}\vert\leq (1/2) \vert z\vert$ then
$$   \vert q(z,\zeta_{1})- q(z,\zeta_{2})\vert \leq C \vert \zeta_{1}-\zeta_{2}\vert (\vert \zeta_{1}\vert + \vert \zeta_{2}\vert)\vert z\vert^{-3/2}. $$

Thus if we write $h=(h_{\bC}, h_{\bR})$ for the components in $\bC\times \bR$ and assume that $\vert h\vert \leq s^{2}$ we have 
$$   V_{0}(x_{0}+h)- V_{0}(x_{0}) (s)=  s^{-1} h_{\bC} + Q(h)(s) $$
where 
$$   \vert Q(h)(s)\vert \leq C \vert h\vert^{2} s^{-3} $$
and $$  \vert Q(h_{1})(s)-Q(h_{2})(s) \vert \leq C \vert h_{1}-h_{2}\vert (\vert h_{1}\vert +\vert h_{2}\vert) s^{-3}.  $$

Our equation (12) becomes
\begin{equation}    L h = Q(h) + W(x_{0}+h)  \end{equation}
where $L$ is the linear operator $Lh= h'-s^{-1} (h_{\bC},0)$. We have an inverse operator
$S$ to $L$ which acts separately on the $\bC$ and $\bR$ components:
$$  S(\rho_{\bC}, \rho_{\bR})=( S_{\bC}(\rho_{\bC}), S_{\bR}(\rho_{\bR}))$$
where  $$ (S_{\bC} \rho_{\bC})(\sigma) = \sigma \int_{0}^{\sigma} s^{-1} \rho_{\bC}(s) ds, $$
and 
$$  (S_{\bR} \rho_{\bR})(\sigma) =  \int_{0}^{\sigma}
\rho_{\bR}(s) ds. $$

We set $h=S(\rho)$ in (13), so the equation becomes the fixed point equation $\rho={\cal F}(\rho)$ where

$$  {\cal F}(\rho) =  Q(S(\rho)) + W(x_{0}+ S\rho). $$

For a given $s_{0}$ we define a weighted norm on functions on $[0,s_{0}]$

$$   \Vert \rho\Vert = \sup s^{-2} \vert \rho(s)\vert . $$

Then $$\vert S(\rho)(s)\vert \leq \frac{1}{2} \Vert \rho\Vert s^{3}. $$
Thus if $\Vert \rho\Vert \leq 1$ (say) then $h=S(\rho)$ satisfies the condition $\vert h\vert \leq s^{2}$, provided that $s_{0}$ is sufficiently small. It is now straightforward to find a fixed point using the contracting mapping theorem and we get uniqueness in the standard way.

\begin{figure}[t]
\centering
\includegraphics[scale=.85]{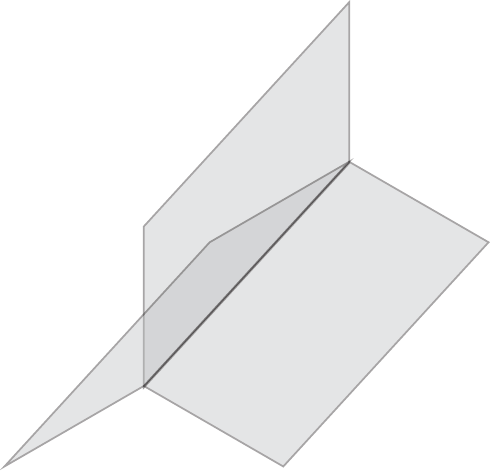}
\caption{Terminal submanifolds meeting at the link.}
\label{fig:terminals}
\end{figure}

For small $t$ we have a similar integral curve $(s^{2},t) + h_{t}(s)$ of the vector field $V$ passing through $(0,t)$. For $x\geq 0$ write
$$ \Phi(x,t)= (x,t)+ h_{t}(\sqrt{x}). $$
Then $\Phi$ is a $C^{1}$ map on a neighbourhood of the origin in the manifold with boundary $\bR^{+}\times \bR\subset \bC\times \bR$ and its image is a $C^{1}$-submanifold $T_{1}$ with boundary
on the axis $z=0$. Each point of $T_{1}$ lies on a gradient line  terminating on the axis. There are two other submanifolds $T_{2}, T_{3}$ obtained in the same way starting with the model solutions $\omega s^{2}, \omega^{2} s^{2}$ where $\omega= e^{2\pi i/3}$. The general picture is much the same as the standard picture for a Morse-Bott function  in three dimensions with a $1$-dimensional critical submanifold and Hessian of index (1,1) in the normal plane. In that case one would have ascending and descending submanifolds near the critical set. In our case the multivalued nature of the function mans that we cannot consistently differentiate between ascending and descending, so we call these {\it terminating submanifolds}. See Figure \ref{fig:terminals}.

Using similar techniques one can show  that the only flow lines (for the vector field corresponding to the vanishing cycle) which terminate on the link $L$ are those that we have constructed above.

\

{\bf Remark} Locally on $L$, there are three terminating submanifolds but it seems possible that these could be permuted as we move around  a component of $L$.

\subsection{Expected dimension and transversality}

The set of gradient cycles has  \lq\lq expected dimension'' $0$. This notion can be made precise in various ways. One approach would be to set up a Fredholm theory in suitable infinite dimensional spaces where the expected dimension appears as the index of a linearised operator. The fact that the dimension is $0$ is essentially  a consequence of the local variational description discussed in Section 5. But since we are studying ordinary differential equations we can also proceed with more elementary methods. For example suppose that we have a gradient cycle formed by a single arc $\gamma$ running from one component $L'$ of $L$ to another component $L''$. Fix a disc $D$ transverse to $\gamma$ at an interior point $q$. Near $L'$ the arc $\gamma$ lies in a $2$-dimensional terminating manifold $T'$. We follow the gradient flow along paths close to $\gamma$ and extend $T'$ until it intersects $D$ in a $1$-dimensional submanifold $S'\subset D$. Similarly we have another submanifold $S''\subset D$ extending a terminating manifold $T''$ near $L''$. By construction $S',S''$ intersect at $q$ and the intersection points $S'\cap S''$ correspond to gradient cycles near $\gamma$. So the meaning of \lq\lq expected dimension'' in this case is that if $S',S''$ are transverse then their intersection has dimension $0$. Of course this is  the same as the usual discussion in Morse-Bott theory of gradient flow lines between critical submanifolds. Similarly, for  a gradient cycle $\Gamma$ which is an orbit we take a transverse disc $D$ and represent nearby gradient cycles as the fixed points of a return map $\phi:D\rightarrow D$, in a standard way.

Now let $\Gamma$ be a gradient cycle for a section $u_{0}$ which is a graph with all vertices of valence $3$.  Let $V$ be the set of vertices $q_{i}$ and $E$ be the set of edges $\gamma_{a}$ and let $I\subset E\times V$ be the incidence set. Choose a disc $D_{a}\subset N$ transverse to the flow line $\gamma_{a}$.  Parametrise $\gamma_{a}(s)$ so that $s=0$ corresponds to the intersection point with $D_{a}$ and suppose that the two end points correspond to $s=l^{+}_{a}, s=-l^{-}_{a}$ for $l^{\pm}_{a}>0$. Let $$\widetilde{D}_{a}= D_{a}\times (-\epsilon,\epsilon)^{2}. $$
For a triple $(z,\eta_{1},\eta_{2})$ in $\widetilde{D}_{a}$ we have two points in $N$ given by following the flow line through $z$ to parameter values $l^{+}_{a}+\eta_{1}$ and $-l^{-}_{a}+\eta_{2}$.  Let $B_{i}$ be a small ball around the point in $N$ corresponding to the vertex $q_{i}$ and for $(a,i)\in I$ let $B_{a,i}$ be a copy of $B_{i}$. The construction above gives a map
$$ f: \prod_{a\in E} \widetilde{D}_{a}\rightarrow \prod_{(a,i)\in A} B_{a,i}. $$
By hypothesis, for each vertex $q_{i}$ there are exactly three pairs $(a,i)$ in $I$ say $(a_{1}, i), (a_{2}, i),(a_{3}, i)$. Let $\Delta_{i}$ be the diagonal in $B_{a_{1}, i}\times B_{a_{2}, i}\times B_{a_{3},i}$,  so we have a submanifold
$$ \Delta^{*}= \prod \Delta_{i}\subset \prod_{I} B_{a,i}. $$
By construction gradient cycles close to $\Gamma$ correspond to points in $\prod \tilde{D}_{a}$ which map by $f$ to $\Delta^{*}$ and we have one such point ${\bf O}$ corresponding to $\Gamma$. Let $n$ be the number of the vertices in $\Gamma$ and $m$ the number of edges. Then  $\Delta^{*}$ has codimension $6n$ in $\prod_{A} B_{a,i}$ and $\prod \tilde{D}_{\alpha}$ has dimension $4m$. But since the graph is trivalent we have $2m=3n$ so $4m=6n$ and if $f$ is transverse to $\Delta^{*}$ at ${\bf O}$ the set of nearby gradient cycles has dimension $0$.

When the transversality condition holds for all gradient cycles of $u_{0}$ we say that the set of gradient cycles is \lq\lq cut out out transversally''.

We expect that for a single section $u_{0}$ the transversality conditions considered above may fail and the set of gradient cycles might have some  more complicated structure but, in the usual way in such theories, we can hope to achieve transversality after suitable generic perturbations. There are many kinds of perturbations we could consider; here we will discuss perturbations of the section $u_{0}$, dropping the maximality condition. (Doing this we lose the \lq\lq Floer-type'', variational, description: there are other kinds of perturbations one could use which retain that.) Thus we let ${\cal H}$ be a space of small perturbations  $u_{0}+\sigma$ with $\sigma$ supported away from the link $L$. (For most purposes we can work with a suitable large finite dimensional space of perturbations.) Then in the setting above our map $f$ extends to
$$  F: {\cal H}\times \prod \tilde{D}_{a}\rightarrow \prod_{A}B_{a,i}, $$ where we use the gradient lines defined using  $u_{0}+\sigma\in {\cal H}$. 

\begin{prop}
The map $F$ is transverse to $\Delta^{*}$ at $(u_{0},{\bf O})$.
\end{prop}

To see this we need to compute the derivative of $F$ with respect
 to $u$. For this we utilize the \lq\lq pertubation theorem'' for flows (see 32.1 in \cite{kn:AR} for example). Let $\Phi^{s}$ be the flow of a smooth vector field $X$ on a manifold, defined near a point $x$ and write $\gamma(s)=\Phi^{s}(x)$. Consider a vector field $Y$ defined on a neighbourhood of $\gamma(s)$ and let $\Phi_{t}^{s}$ be the flow of  $X+tY$. Then
\begin{equation}  \frac{d}{dt} \Phi^{s}_{t}(x)\vert_{t=0}= \int_{0}^{s} (\Phi^{s})_{*}(Y_{\gamma(s-r)}) dr. \end{equation}
More generally the same formula holds for a flow defined by a $1$-parameter family of vector fields $X_{t}$ with $t$-derivative $Y$ at $t=0$.

Let $\gamma(s)=\gamma_{a}(s)$ be the flow line corresponding to an edge of $\Gamma$, as considered above, with $s=0$ at the intersection point $0$ with the transverse disc $D_{a}$. Write $l^{+}=l^{+}_{a}$ so that $q_{+}=\gamma(l^{+})$ is a vertex  of $\Gamma$. We consider first perturbations $u_{0}+ t\sigma$ where $\sigma$ is supported in a small ball  around $\gamma(l^{+}/2)$ (say). We get perturbed gradient lines
$\gamma_{t}(s)$ passing through the same point $0$. For $s\leq 0$ we have $\gamma_{t}(s)= \gamma(s)$ but for $s>0$ they are different. Let $V$ be the derivative
$$  V= \frac{d}{dt} \gamma_{t}(l^{+}).$$
So $V$ is a vector in  the tangent space of $N$ at $q_{+}$. 

\begin{lem}
We can choose a pair of such pertubations $\sigma_{1}, \sigma_{2}$ giving vectors $V_{1}, V_{2}$ in $TN_{p_{+}}$ such that $V_{1}, V_{2}, \gamma'(l^{+})$ form a basis of  $TN_{p_{+}}$.
\end{lem}

Assuming this Lemma we deduce Proposition 3 as follows. Focusing on the vertex $q_{+}$ and the edge $\gamma_{a}$, we have a ball $B$ around $q_{+}$ and a map $$F^{+}_{a}: {\cal H}\times D_{a}\times (-\epsilon,\epsilon)\rightarrow B,  $$ which gives one component of the map $F$.
The derivative of $F^{+}_{a}$ with respect to $\eta\in (-\epsilon,\epsilon)$ is given by $\gamma'(l^{+})$ and the derivative with respect to the two variations provided by the Lemma gives $V_{1}, V_{2}$. The derivatives of all other components of the map $F$ with respect to these variations vanish so taking $4 m$ similar variations (two for each pair $(a,i)\in A$) we see that the derivative of $F$ is surjective so certainly $F$ is transverse to $\Delta^{*}$.

We now prove the Lemma. Recall that in a local description, the gradient vector fields we are considering are given by the projection of a fixed vector $\alpha$ in $\bR^{3,19}$ to the tangent space of the image of the map $u$ from a neighbourhood $B$ in $N$ to $\bR^{3,19}$. We take  $u=u_{0}+ t \sigma$ and choose $\sigma$ to be a normal variation so that for each $x\in B$ the value $\sigma(x)$ is orthogonal to the image of $du$ in $\bR^{3,19}$. One gets the following formula for the $t$-derivative $Y$ of the gradient vector field $X$:
\begin{equation} Y = {\rm grad} \langle \sigma, \alpha\rangle - \sigma S X \end{equation} where $S$ is the second fundamental form of the image submanifold. In our application using the formula (14) we only need the value of $Y$ on the curve $\gamma$ so we choose $\sigma$ to vanish on $\gamma$ and this means that the second term in the formula (15) vanishes. The derivative
$   (\Phi^{s})_{*}$ appearing in the formula (14) preserves the tangent vectors to $\gamma$. so for our purposes we only need to consider the component of $Y$ normal to $\gamma$. It is straightforward to construct a $\sigma$ of this kind realising {\it any} normal vector field  and then to use (14) to produce the desired $\sigma_{1}, \sigma_{2}$.

The whole discussion above can be adapted without difficulty to general gradient cycles and using standard techniques we obtain:

\begin{cor}

There is a residual set ${\cal P}$ of positive sections such that for any $u\in {\cal P}$ the set of gradient cycles is a $0$-manifold cut out transversally.  \end{cor}

The results of this section can be extended to many situations when the cycles are allowed self-intersections although there are some difficulties in proving a completely general result.

\subsection{Crossing}
We have now reached the main point of this article, which is to study phenomena for gradient cycles in $1$-parameter families analogous to those we outlined in Section 2 for associative submanifolds. In the remainder of  Section 6 we will be content to make whatever transversality assumptions are relevant, without analysing the exact meaning of these. 
In this subsection we study what happens when, in a generic $1$-parameter family,  two gradient graphs cross in the complement of the link.

Let $\Gamma\subset N$ be a gradient graph for a section $u_{0}$ and $q$ a point of $\Gamma$ which is not a vertex, so $q$ lies in the interior of some gradient flowline $\gamma$, an integral curve of a local vector field $V$. We introduce a notion of \lq\lq cut'' gradient graphs. Let $D$ be a small disc in $N$, centred at $q$ and transverse to $\gamma$. For each $z$ in $T$ there is is gradient flowline $\gamma_{z}$ through $z$.  The intersection of $\Gamma\setminus \{q\}$ with a small ball centred at $q$ has two components.  By a cut gradient graph we mean a small deformation of $\Gamma$ away from $q$ but near $q$ we allow these two components to deform to  possibly different gradient curves $\gamma_{z}, \gamma_{z'}$. We have a moduli space $M$ of such cut gradient graphs which maps to $D\times D$. Our dimension analysis above shows that  $M$ is a $2$-manifold and we suppose that the image of  $M$ in $D\times D$ intersects the diagonal transversally at $(q,q)$. So we can assume that $M$ is embedded in $D\times D$. 
  
Next suppose that $\Gamma_{1},\Gamma_{2}$ are two such gradient graphs (for the same section $u_{0}$) which intersect at a point $q$ (not a vertex of either graph). So $q$ lies on $\gamma_{1}$ and $\gamma_{2}$ which are integral curves of vector fields $V_{1}, V_{2}$ which we assume are linearly independent at $q$. In our situation we want to assume that the $V_{i}$ are defined by $-2$ classes $\alpha_{i}$ with $\alpha_{1}.\alpha_{2}=1$ but this condition will not play a fundamental role in the discussion below.
We choose transverse discs $D_{1}, D_{2}$ as above and we have moduli spaces $M_{1}, M_{2}$ of cut graphs. Let $I\subset D_{1}\times D_{2}$ be the set corresponding to  pairs of intersecting flow lines. It is clear that $I$ is a $3$-dimensional submanifold of $D_{1}\times D_{2}$ and the intersection point defines a map $f$ from $I$ to $N$. In the $8$-manifold $$D_{1}\times D_{2}\times D_{1}\times D_{2}= D_{1}\times D_{1}\times D_{2}\times D_{2}$$ we have two subsets $I\times I$ and $M_{1}\times M_{2}$. Let $J$ be their intersection. We assume that this intersection is transverse, so $J$ is a $2$-manifold and $F=(f\times f)$ gives a map $F:J\rightarrow N\times N$. By our hypotheses the point $(q,q)$ lies in the image of $F$ : it is the image of the point $Q= (q,q,q,q)$ in $J$.

We now introduce a $1$-parameter family of sections $u_{t}$ for $\vert t\vert<\epsilon$ so $\Gamma_{1}, \Gamma_{2}$ deform in families $\Gamma_{i}^{t}$. We can deform all the constructions above in the family, so for each $t$ we have a $2$-manifold $J_{t}$ and a map $F_{t}: J_{t}\rightarrow N\times N$. The difference is that for $t\neq 0$ we do not expect that the image of $F_{t}$ meets the diagonal in $N\times N$.  We can choose a family of co-ordinate charts $\psi_{t}: B^3\rightarrow N$ which linearise the vector fields $V_{1}^{t}+ V_{2}^{t}$. In other words, $V_{1}^{t}+V_{2}^{t}$ is the image by $d\psi_{t}$ of a constant vector field $n$ on the unit ball $B^3$ in $\bR^{3}$. Let $G_{t}:J_{t}\rightarrow \bR^{3}$ be the composite of $(\psi_{t}^{-1}\times \psi_{t}^{-1})\circ F_{t}$ with the difference map $(x,y)\mapsto x-y$ from $\bR^{3}\times \bR^{3}$ to $\bR^{3}$.

\begin{figure}[t]
\centering
\includegraphics[scale=.85]{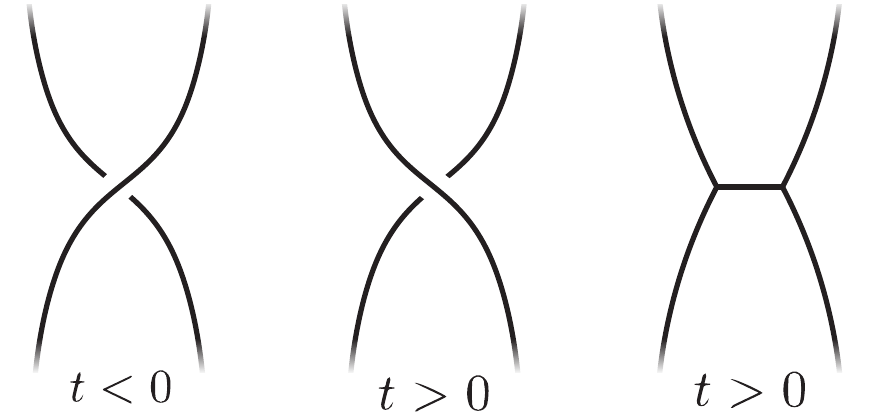}
\caption{An adiabatic analogue of Joyce-Nordst\"{o}m crossing}
\label{fig:crossing}
\end{figure}

Let ${\cal J}$ be the $3$-manifold formed by the family of $2$-manifolds $J_{t}$,  so we have a map $\pi:{\cal J}\rightarrow (-\epsilon,\epsilon)$ with fibres the $J_{t}$ and a map ${\cal G}: {\cal J}\rightarrow \bR^{3}$, equal to $G_{t}$ on the fibres.  By construction we have a point  $Q$ in $J_{0}\subset {\cal J}$ with ${\cal G}(Q)=0$. We make the transversality assumptions:
\begin{enumerate}\item  $Q$ is a regular point of ${\cal J}$,
\item the vector $n$ does not lie in the image of $dG_{0}$ at $Q$.
\end{enumerate}
Let $Z\subset \bR^{3}$ be the ray generated by $n$. These transversality assumptions imply that the pre-image ${\cal G}^{-1}(Z)$ is a $1$-manifold with boundary embedded in ${\cal J}$ and it intersects the fibres $J_{t}$ either for small positive $t$ or for small negative $t$ (but not both). Without loss of generality,  suppose that the intersection occurs for small positive $t$. Then, by construction, for such $t$ there is a gradient graph $\Gamma^{t}$ of the kind indicated in the right of Figure \ref{fig:crossing}. In other words, in the family a new gradient graph is \lq\lq born'' at $t=0$.

The above discussion gives an adiabatic analogue of the Joyce-Nordstr\"om crossing for associative submanifolds outlined in Section 2. One can see that the topology matches up. Recall that for each gradient graph $\Gamma$ we have defined, at the topological level, a submanifold $P(\Gamma)\subset M^{7}$.  In the situation above, the fact that $\alpha_{1}.\alpha_{2}=1$ implies that $P_{i}= P(\Gamma_{i})$ intersect in one point and one can check that the submanifold $P(\Gamma^{t})$ corresponding to the graph $\Gamma^{t}$ is the connected sum $P_{1}\sharp P_{2}$.

\subsection{Splitting of flowlines and surgery triples}

In this subsection we discuss an adiabatic analogue of the \lq\lq surgery triple'' outlined in Section 2. The essential phenomenon arises for general pairs of vector fields on a manifold. Let $V_{1}, V_{2}$ be non-zero vector fields on the open unit ball $B^{3}$ and $\Phi^{s}_{1}, \Phi^{s}_{2}$ be the local flows they generate. Fix $\delta>0$ so that $\Phi_{i}^{s}$ are defined with values in $B^{3}$ for $s\leq 2 \delta$ on the half-sized ball $\frac{1}{2}B^{3}$. Define 
$$   C= \{ x\in \frac{1}{2} B^{3}: V_{1}(x)= \lambda V_{2}(x)\  {\rm for \  some}\ \lambda>0\}, $$
and $$ S= \{(x, x')\in \frac{1}{2} B^{3}\times B^{3} : \Phi_{1}^{s_{1}}(x)=x', \Phi_{1}^{s_{2}}(x)=x' \ {\rm for
\  some} \ 0<s_{1}, s_{2}<\delta\}. $$
So we have an inclusion map $F:S\rightarrow  B^{3}\times B^{3}$ and we can consider $C$ as contained in the diagonal  $B^{3}\subset B^{3}\times B^{3}$.

The set $C$ corresponds to the intersection of the two sections of the unit sphere bundle defined by $V_{i}/\vert V_{i}\vert$. We assume that this intersection is transverse, so that $C$ is a $1$-dimensional submanifold of $B^{3}$. 
\begin{prop}
Under this transversality assumption, $C$ lies in the closure $\overline{S}$ of $S$ and near $C$ the closure $\overline{S}$ has the structure of a $2$-manifold with boundary $C$.  
\end{prop}

To prove the Proposition, let $q$ be a point of $C$ and choose local co-ordinates $(x_{0}, x_{1}, x_{2})$ centred on $q$ so that $V_{1}=\frac{\partial}{\partial x_{0}}$. The statement of the Proposition is unchanged if we multiply $V_{2}$ by a positive function so we can suppose that 
$$  V_{2}= \frac{\partial}{\partial x_{0}}+\xi_{1}  \frac{\partial}{\partial x_{1}}+ \xi_{2}  \frac{\partial}{\partial x_{2}}, $$
where $\xi_{i}$ are functions of $x_{0}, x_{1}, x_{2}$ vanishing at the origin.
The transversality condition is that the the $2\times 2$ matrix $\frac{\partial \xi_{i}}{\partial x_{j}}$ ($i,j=1,2$) is invertible. The submanifold $C$ is the common zero set of $\xi_{1}, \xi_{2}$.  The local flow $\Phi_{1}^{s}$ of $V_{1}$ is just translation in the $x_{0}$ factor and we have  
$$  \Phi_{2}^{s}(0,{\bf x}) = (s, \phi_{s}({\bf x}))$$ for a family of diffeomeorphsims $\phi_{s}$ of neighbourhoods of $0$ in  $\bR^{2}$. If $a<b$ and ${\bf x}\in  \bR^{2}$ is a fixed point of the diffeomorphism $\phi_{b}\circ \phi_{a}^{-1}$ then the pair $\left((a,{\bf x}),(b, {\bf x})\right)$ lies in $S$ and all points of $S$ arise in this way. The Proposition is now a consequence of the following simple fact. Suppose that $\psi_{s}$ is a family of diffeomorphisms of a neighbourhood of $0$ in $\bR^{2}$ with $\psi_{0}$ equal to the identity. Let $\eta$ be the vector field given by the $s$ derivative of $\psi_{s}$ at $s=0$. If $\eta$ has a transverse zero at the origin then for small $s$ there is a unique fixed point of $\psi_{s}$ close to the origin. Applying this to the $\phi_{b}\phi_{a}^{-1}$ we get a local parametrisation of $S$ by pairs $(a,b)$ with $a<b$ which extends to the diagonal $a=b$, mapped to $C$.

\begin{figure}[t]
\centering
\includegraphics[scale=.8]{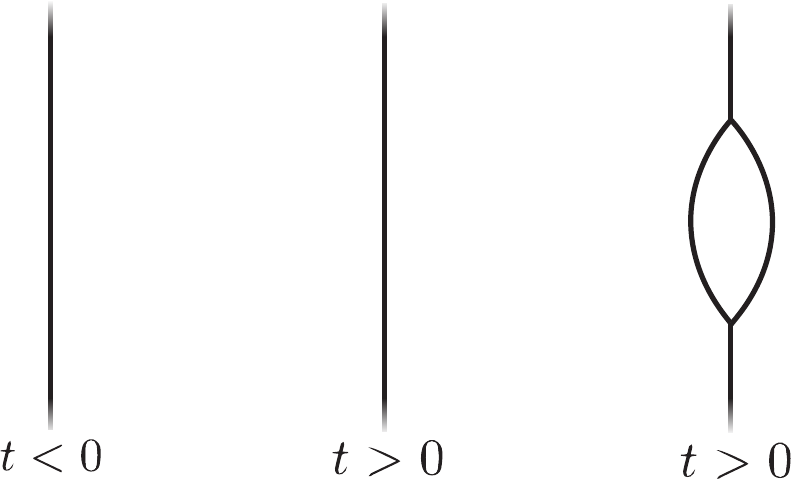}
\caption{An adiabatic analogue of a surgery triple}
\label{fig:surgery}
\end{figure}

Let $V$ be the vector field $V_{1}+V_{2}$. Let $\gamma$ be the flow line of $V$ through the origin and choose a transversal disc $D$ as before, parametrising nearby  flow lines $\gamma_{z}$ of $V$. The quotient by the flow defines a map $p:B^{3}\rightarrow D$ so we get a map
$$G= (p\times p)\circ F :S\rightarrow D\times D. $$
This extends to $\overline{S}$, mapping $C$ to the diagonal in $D\times D$.
Let $M$ be a $2$-dimensional submanifold of $D\times D$. If $(z,z')\in D\times D$ is an intersection point of $M$ and $S$ we have a configuration of flow lines of the kind seen on the right of Figure \ref{fig:surgery}. This consists of a segment of a flow line $\gamma_{z}$ for $V$ ending at a point $x$, a pair of   flow lines of $V_{1}, V_{2}$ from $x$ to another point $x'$ and a segment of the flow line $\gamma_{z'}$ starting from $x'$.

Suppose now that the origin in $B^{3}$ lies in the submanifold $C$. Let $V_{1}^{t}, V_{2}^{t}$ be a $1$-parameter family of deformations, so for each $t$ we have  submanifolds $C_{t}, S_{t}$. Let ${\cal C}, {\cal S}$ be the corresponding sets in $B^{3} \times B^{3}\times (-\epsilon,\epsilon)$ so we have a map
$$  {\cal G}:\overline{{\cal S}}\rightarrow D\times D\times (-\epsilon,\epsilon), $$ which maps the boundary ${\cal C}$ to the diagonal times $(-\epsilon,\epsilon)$.
Similarly, a family $M_{t}\subset D\times D$ defines a $3$-dimensional submanifold
${\cal M}\subset D\times D\times (-\epsilon,\epsilon)$.  By construction the image of ${\cal G}$ intersects ${\cal M}$ at the point ${\bf 0}=(0,0,0)$. We make the transversality assumption that this  intersection is transverse; then the intersection is a $1$-dimensional manifold $I$ with boundary ${\bf 0}$   and the derivative of the projection map $I\rightarrow (-\epsilon,\epsilon)$ at ${\bf 0}$ is non-zero. Depending on the sign of this derivative there is  a configuration of the kind described above for the vector fields $V_{i}^{t}$ and submanifolds $M_{t}$ {\it either} for small positive $t$ {\it or} for small negative $t$, but not both. Without loss of generality we suppose the first case.

 In our situation we take $B^{3}$ to be a co-ordinate chart centred at a point $q$ in $N\setminus L$ and $V$ to be the gradient vector field defined by a  $-2$ class $\alpha$. We suppose that $q$ lies on a segment $\gamma$ of a cycle $\Gamma$ which satisfies all the conditions to be a gradient cycle except that the irreducibility condition fails at $q$. We let $M$ be the submanifold parametrising cut gradient cycles, as before. The central assumption is that $\alpha=\alpha_{2}+\alpha_{3}$ for $-2$ classes $\alpha_{i}$ giving gradient vector fields $V_{i}$ with $V_{1}=\lambda V_{2}$ at $q$. Then we see that as we vary our section $u$ in a $1$-parameter family $u_{t}$, satisfying the transversality hypotheses, there will be a family of gradient cycles $\tilde{\Gamma}_{t}$ \lq\lq born'' at $t>0$. See Figure \ref{fig:surgery}.

\

We interpret this as the adiabatic analogue of a surgery triple obtained by the three Harvey-Lawson smoothings of a  cone singularity. We get evidence for this from the topology of the situation. The condition that $V_{1}=\lambda V_{2}$ at $q$ means that we are in the exceptional case where the complex curve $\Sigma\subset X$ is reducible, a union of two embedded spheres $\Sigma_{1}, \Sigma_{2}$ intersecting in a point in $X$. For a nearby point $q'$ in $\Gamma$ the corresponding curve $\Sigma'\subset X$ is smooth but \lq\lq close'' to singular,  with a small circle $\lambda\subset \Sigma'$ (the vanishing cycle) which shrinks to a point as $q'$ moves towards $q$. Now consider a small  deformation of $\Gamma$ to $\Gamma_{t}$ for $t<0$. We have a $3$-manifold $P_{-}= P(\Gamma_{t})$ in $M^{7}$ and we can deform $\lambda$ slightly to a circle $\lambda_{-}$ in $P_{-}$. We also have a $3$-manifold $\tilde{P}$ defined by $\tilde{\Gamma}_{t}$ for $t>0$. We leave the reader to check that $\tilde{P}$ is obtained from $P_{-}$ by $0$-surgery on $\lambda_{-}$. The third $3$-manifold in the picture is given by $P_{+}= P(\Gamma_{t})$ for $t>0$, which contains a small circle $\lambda_{+}$. In fact $P_{+}$ is diffeomorphic to $P_{-}$ but we can consider a more refined notion. Let $U_{\pm}$ be small tubular neighbourhoods of $\lambda_{\pm}$ in $P_{\pm}$. There is a natural diffeomorphism $F: P_{-}\setminus U_{-}\rightarrow P_{+}\setminus U_{+}$, well-defined up  to a small isotopy. Then, relative to $F$, the $3$-manifolds $P_{+}, P_{-}$ are different: i.e. we cannot extend $F$ to a diffeomorphism from $P_{-}$ to  $P_{+}$. Relative to $F$,  and with suitable orientations, the manifold $P_{+}$ is obtained from $P_{-}$ by $+1$ surgery on $\lambda_{-}$. This follows from the fact that if we have a standard family of complex structures $I_{\tau}$ on $X$ parametrised by $\tau$ in the unit disc in $\bC$ such that the $-2$ class $\alpha$ is represented by a smooth curve for $\tau\neq 0$ and the singular curve $\Sigma_{1}\cup\Sigma_{2}$ for $\tau=0$ then  the monodromy of the family of curves is given by the Dehn twist in the vanishing cycle $\lambda$. This Dehn twist is trivial in the mapping class group of $S^{2}$ but is non-trivial in the appropriate relative mapping class group.

\subsection{Other transitions}

In this subsection we outline (without full proofs) two phenomena involving gradient flow lines near the link $L$. The first is another version of  \lq\lq crossing''. Suppose that $\Gamma_{1}, \Gamma_{2}$ are two gradient graphs for a section $u_{0}$ each of which has a vertex of valence $1$ at the same point $p$ of the link $L$. As usual we consider a $1$-parameter family of sections $u_{t}$ and suppose that $\Gamma_{i}$ deform in families $\Gamma_{i}^{t}$. So we have vertices $p_{1}(t), p_{2}(t)$ say on $L$ which coincide at $t=0$. We make the transversality assumption that the map $t\mapsto (p_{1}(t), p_{2}(t))$ is transverse to the diagonal at $t=0$. Then one can show that  there is another family of gradient graphs $\left(\Gamma_{1}\sharp\Gamma_{2}\right)^{t}$ defined either for small positive $t$ or for small negative $t$, but not both. The gradient graph is obtained by deforming the union of two flow lines terminating on $L$ into a single flow line. The basic model is given by
the gradient curves of the function ${\rm Re}( z^{3/2})$ with respect to the Euclidean metric, as in Subsection 6.1. Then, up to  parametrisation, there are gradient curves 
$$  \gamma_{\epsilon}(s)= (s+i \epsilon)^{2/3} , $$
which converge in an obvious sense as $\epsilon\rightarrow 0$ to the union of two line segments. This can be seen in any of the three sectors in Figure \ref{fig:vector-field}. It is easy to see that $P((\Gamma_{1}\sharp\Gamma_{2})^{t})$ is topologically the connected sum of $P(\Gamma_{1})$ and $P(\Gamma_{2})$. 

\

For the second phenomenon we consider a gradient cycle $\Gamma_{-\epsilon}$ which contains a segment $\gamma$ that passes close to a component $L_{0}$ of the link $L$.  We have two $-2$ classes $\alpha,\delta$ where $\alpha$ is the label of $\gamma$ and $\delta$ is the vanishing cycle and we suppose that $\alpha.\delta=1$. The multi-valued setting means that the gradient vector field $V$ associated to $\alpha$ is not well-defined on a neighbourhood of $L_{0}$, even up to sign. We can define $V$ on a  cut neighbourhood but there will be a jump by the addition of the vector field $V_{\delta}$ associated to $\delta$ across the cut. Since $V_{\delta}$ vanishes on $L_{0}$ the value of $V$ at points of $L_{0}$ is well-defined and is not zero. 

The basic phenomenon can be seen in the model case where we work on $\bC$ with the gradient vector field $V$ of the multivalued function ${\rm Re}(z^{3/2}+ iz) $ which we interpret by making a cut along the positive real axis. The jump across the cut is by $V_{\delta}= \frac{3}{2} z^{1/2}$. There is a family of flow lines $\gamma_{t}$ for $t<0$ given up to parametrisation by ${\rm Im}(z^{3/2}+ i z)=t$. These do not meet the cut. The function ${\rm Im}(z^{3/2}+ iz )$ extends continuously across the cut but its derivative does not. If we attempt to extend the definition to $t>0$ the same equation defines a set which meets the cut at the point $p_{t}=t$. But $p_{t}$ also lies on a flow line of $V_{\delta}$ through the origin. Instead of a single flow line we get a configuration of three flow lines meeting at the point $p_{t}$. If we change our point of view and define the multivalued function by making a cut along the negative real axis (say) then the flow lines are associated to the three vector fields $V, V_{\delta}, V+V_{\delta}$.

\begin{figure}[t]
\centering
\includegraphics[scale=.8]{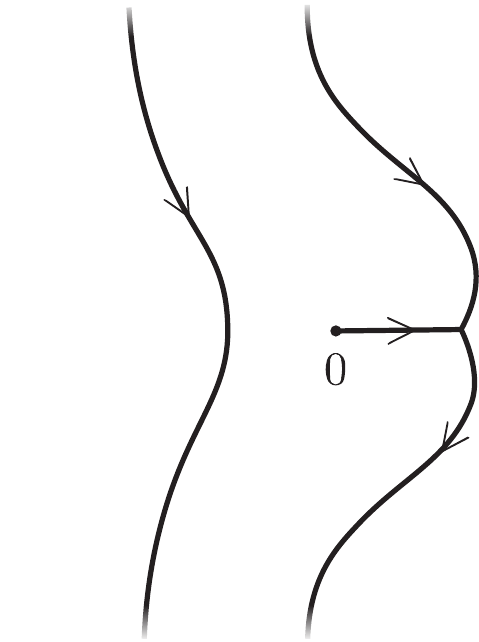}
\caption{A gradiant cycle passes through the link $L$ in a cross-section, where $L$ corresponds to the origin in $\bC$.}
\label{fig:cut-plane}
\end{figure}

In a similar fashion, for a $1$-parameter family of sections $u_{t}$ we may see a gradient cycle $\Gamma_{-\epsilon}$ for $t=-\epsilon$ which deforms
in a family $\Gamma_{t}$ where $\Gamma_{0}$ meets $L_{0}$. This family $\Gamma_{t}$ can be extended to $t>0$ but with a graph of a different topological type, inserting a vertex $p_{t}$ of valence $3$ and a vertex of valence $1$ on $L_{0}$. The labels of the three edges meeting at $p_{t}$, with outgoing orientations,  are $-\alpha, \alpha+\delta, '-\delta$. The essential point here is that the monodromy of the flat bundle around $L_{0}$ takes $\alpha$ to $\alpha+\delta$. 

In this case we do not change the \lq\lq count'' of gradient cycles and the $3$-manifolds $P(\Gamma_{t})$ are homeomorphic for $t$ positive and negative.



\end{document}